\documentclass[11pt]{amsart}
\usepackage{amsthm,amssymb,url,hyperref,fullpage,graphicx}
\usepackage[usenames]{color}

\theoremstyle{plain}
\newtheorem{theorem}{Theorem}[section]

\theoremstyle{definition}
\newtheorem{example}[theorem]{Example}

\def\lmlt#1{\mathrm{LMlt}(#1)}

\def\dis#1{\mathrm{Dis}(#1)}
\def\Z{\mathbb Z}

\begin{document}

\title{The origins of involutory quandles}
\author{David Stanovsk\'y}
\address{Department of Algebra, Faculty of Mathematics and Physics, Charles University, Prague, Czech Republic}
\email{stanovsk@karlin.mff.cuni.cz}
\date{Originally written in 2004. Thouroughly revised on \today}

\begin{abstract}
We present an overview of some older papers on involutory quandles, mostly from the times before the term ``quandle" was born. It is meant as a reference guide. It is not meant (yet) as an expository article explaining what the involutory quandles are and what they are good for.
\end{abstract}

\maketitle 

{\bf Note.}
{\it This is an unfinished work and it is very likely incomplete (missing references) or inaccurate (incorrect or unfair description of the results). At the moment, I think it is valuable to post such a reference guide, but it is not worth the effort to polish it to the state of a rigorous mathematical publication.}

\section{Introduction}

An \emph{involutory quandle} is an algebraic structure $(Q,\cdot)$ satisfying the following identities:
\begin{align}
x(xy)&=y,\tag{involutory, aka left symmetric}\\
x(yz)&=(xy)(xz),\tag{left distributivity}\\
xx&=x. \tag{idempotence}
\end{align}
Various other names have been used in literature: \emph{kei} by Takasaki (and ocassionally later) \cite{NP,Ta}, \emph{symmetric sets} by Nobusawa and his followers \cite{No1}, \emph{symmetric groupoids} by Pierce \cite{Pi1}, \emph{left symmetric left distributive idempotent groupoids} by the Prague algebraic group \cite{Sta-phd}.

An \emph{involution} (or involutory mapping) is a permutation $f$ such that $f^2=id$.
Given a binary algebraic structure $(Q,\cdot)$, consider the left translations $L_a:Q\to Q$, $x\mapsto ax$. Then $(Q,\cdot)$ is an involutory quandle if and only if it is idempotent and all left translations are involutory automorphisms. 
If we only assumed that left translations were automorphisms (regardless their order), we would have obtained the class of \emph{quandles}. Many results presented here have a generalization into the class of all quandles.

Our goal is not a compact exposition of the theory of involutory quandles. Rather, we want to tell the story of how the investigations of involutory quandles have been conducted. Therefore, the primary ordering of the material is time. In Section 2, we present three main motivations to study involutory quandles. In Section 3, we overview most(?) of the papers published on involutory quandles in 20th century. In Section 4, we briefly mention how to study connected involutory quandles using a homogeneous representation and the canonical correspondence to transitive groups developed in \cite{HSV}, and show some enumeration results.

The text is nearly self-contained, although familiarity with the basics of quandle theory, as described e.g. in \cite[Sections 2,3]{HSV} or in \cite[Section 1]{AG}, would improve reader's understanding.
Let $\lmlt Q=\langle L_a:a\in Q\rangle$ and $\dis Q=\langle L_aL_b:a,b\in Q\rangle$ be the (left) \emph{multiplication group} and the \emph{displacement group} of a quandle $Q$, respectively (aka inner group and transvection group, resp.) .
A quandle $Q$ is called \emph{connected} if $\lmlt Q$ (or, equivalently, $\dis Q$) acts transitively on $Q$. It is called \emph{faithful} (aka effective) if all left translations are pairwise different. It is called \emph{latin} (aka a quasigroup) if all right translations are permutations, too. It is called \emph{medial} (aka abelian or entropic) if it satisfies the identity $(xy)(uv)=(xu)(yv)$.

\section{Motivation}

\subsection{Symmetric spaces: Abstraction of reflection}

In a euclidean space, let $a*b$ denote the reflection of $b$ over $a$. The resulting algebraic structure is an involutory quandle. This observation, and the resulting abstraction of the notion of a reflection, can be attributed to Takasaki and his remote 1942 work \cite{Ta}, but the real advances have been made by Loos and others two decades later \cite{Lo}. 
His \emph{symmetric spaces} can be defined as differentiable manifolds equipped with a smooth involutory quandle operation such that, for any point $a$, there is a neighborhood $U$ of $a$ such that $a$ is the only fixed point of the left translation $L_a$ on $U$.

\begin{example}
On any manifold with a metric, define a symmetric space by $x*y$ to be the image of $y$ by the symmetry through $x$.
For instance, on the real line $x*y= 2x-y$, on the $n$-sphere $x*y=2\left<x,y\right>x-y$. The 3-sphere gives a reason why the fixed point condition is defined only locally: the point symmetric to the south pole through the north pole is the south pole, hence the uniqueness condition fails. 
Another example of a symmetric space is the set of all $n\times n$ positive definite symmetric matrices over the reals with the operation $A\circ B=AB^{-1}A$. Grassman manifolds and Jordan algebras are also symmetric spaces. Any involutory quandle can be considered as a discrete symmetric space. 
\end{example}

An interesting perspective on symmetric spaces is explained in \cite{Kik-motivation}.
Another class of geometric examples of involutory quandles, releated to homotopy classes on spheres, is discussed in \cite{En2}.

\subsection{Cores: Isotopy invariants of loops}

Let $(G,\cdot)$ be a group, or, more generally, a Bol loop (i.e. a "non-associative group" where associativity is replaced by a weaker law $x(y(xz))=(x(yx))z$). The binary algebra $(G,*)$ with \[a*b=a(b^{-1}a)\] is an involutory quandle, called the \emph{core} of $(G,\cdot)$.
Cores were introduced by Bruck who studied isotopy invariants of loops. He proved that isotopic Moufang loops have isomorphic cores \cite{Br}. Cores were later picked up by Belousov and others to construct some of the first examples of latin quandles, see e.g. \cite[Chapter IX]{Bel-book}.

Cores were also investigated on its own.
The equational theory of involutory quandles is that of the group cores \cite{Pi1}. 
The core of a group is 
\begin{itemize}
	\item right distributive iff the group satisfies the identity $xy^2x=yx^2y$ \cite{Bel-cores}.
	\item medial iff the group is 2-nilpotent (an unpublished result by Roszkowska).
	\item faithful iff the center of the group does not contain an involution \cite{Um}. 
	\item latin iff the group is uniquely 2-divisible (i.e. the mapping $x\mapsto x^2$ is a permutation) \cite{Bel-book,KNN}.
\end{itemize}
If the center of $G$ is trivial, one can state stronger results \cite{Um}: 
the displacement group can be embedded into the product of $L\times R$, where $L,R$ are the groups generated by the left and right translations of $G$, respectively, it consists of all $L_bL_aR_{b'}R_a$ such that $a\in G$, $b,b'\in G'$ (here $G'$ denotes the derived subgroup), and if $N$ is the orbit of $\dis Q$ containing the unit, then $N$ is a normal subgroup of $G$ and $G/N$ is an elementary abelian subgroup of exponent 2.
Consequently, if $G=G'$, then $\dis Q$ is isomorphic to $L(G)\times R(G)$.

More about cores of Bol loops can be found in \cite[Chapter IX]{Bel-book} or \cite[Section III.8]{Sta-phd}.
Some variations on the core construction are discussed in \cite{Sta-gop}.

\subsection{Knot quandles: Ambient isotopy invariants of knots}

The following construction is due to Joyce \cite{Jo1} and Matveev \cite{Ma}. 
Let $K$ be a tame knot and consider its diagram. It divides the knot to arcs (segments from one underpass through some overpasses to the next underpass), let $A$ denote the set of arcs. The \emph{involutory knot quandle} is the involutory quandle generated by $A$ subject to the relations $ab=c$ for every crossing with the bridge $a$ and the underpassing arcs $b,c$ (note that the position of $b,c$ plays no role since $ab=c$ iff $ac=b$). The \emph{knot quandle} is a richer structure, taking into account an orientation: the definition is analogous, requiring additionally that $b$ is the right underpass and $c$ is the left underpass.

Joyce and Matveev, independently, proved that the knot quandle is a complete invariant with respect to ambient isotopy (up to mirroring); in particular, the knot quandle does not depend on the diagram.
The involutory knot quandle is an incomplete invariant, nevertheless Winker \cite{Win} proved that it distinguishes the unknot (i.e. non-trivial knots have non-trivial involutory quandles).

\begin{center}
\includegraphics[scale=0.25]{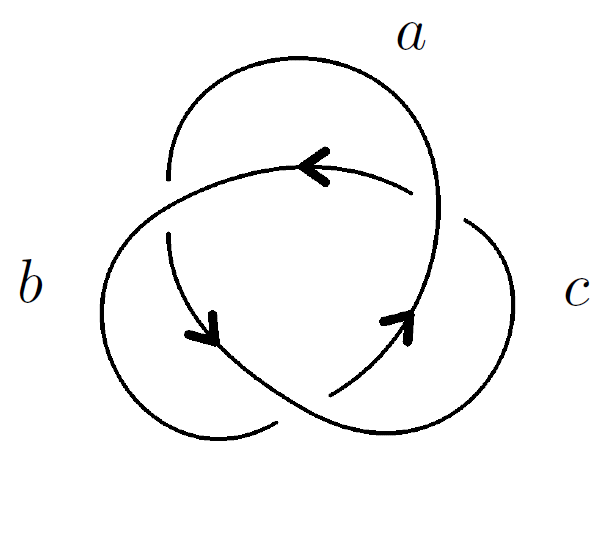}
\hskip3cm
\includegraphics[scale=0.33]{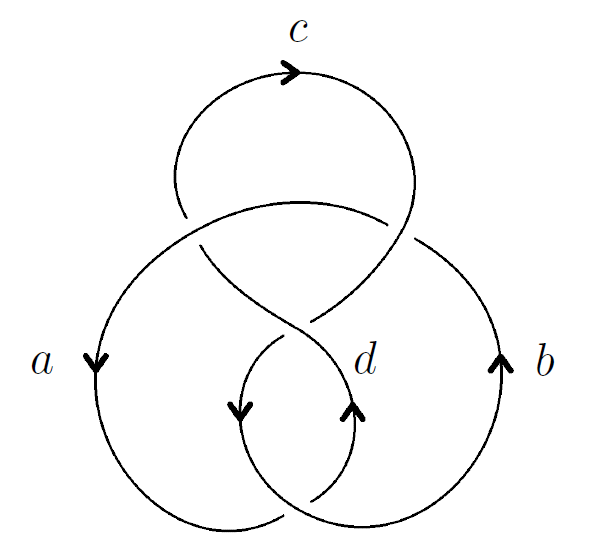}
\end{center} 

\begin{example}
The trefoil quandle is presented by
\[\left<a,b,c;\ ac=b,\ ba=c,\ cb=a\right>\]
and the figure-eight quandle is presented by
\[\left<a,b,c,d;\ ac=d,\ ba=d,\ ca=b,\ dc=b\right>.\]
It is easy to calculate that the involutory quandle of the trefoil is the core of the group $\Z_3$,
while for the figure-eight it is the core of $\Z_5$.
\end{example}

\section{Algebraic results}

\subsection{Nobusawa's project}

Inspired by Loos, Nobusawa started a purely algebraic study of involutory quandles in the 1970s. He used the name \emph{symmetric set} with the following definition: it is a mapping $S$ from a set $A$ into the set of involutions over $A$, sending an element $a$ onto an involution $S_a$ such that $S_a(a)=a$ and
$S_{S_a(b)}=S_aS_bS_a$ for all $a,b\in A$. Then $a*b=S_a(b)$ is an involutory quandle operation on $A$, and every symmetric set results from an involutory quandle as the mapping $a\mapsto L_a$.

In his first article \cite{No1}, Nobusawa formulated some of the fundamentals of the structure theory of involutory quandles. He demonstrated that the displacement group reflects a lot of properties of the quandle. For example, an involutory quandle $Q$ is medial if and only if the group $\dis Q$ is abelian. If $Q$ is faithful, then it is medial if and only if $\dis Q=\{L_eL_a:\, a\in Q\}$ for a fixed $e\in Q$, and in such a case, $Q$ is latin (called homogeneous there). 
Fix $e\in Q$ and define the power of an element $a$ by \[a^k=\underbrace{e(a(e(a\dots)))}_{k+1}.\]
In this way, \emph{cycles} are defined (in other terms, cycles are homomorphic images of the core of the group of integers in $Q$), and their properties investigated (a somewhat clearer account on cycles can be found in \cite[Section III.5]{Sta-phd}).
A finite involutory quandle is latin if and only if all cycles have odd length.
Latin involutory quandles with prime and prime square elements are cores of abelian groups, but there is a latin involutory quandle with 27 elements which is not the core of a group. The last page contains a list of all involutory quandles with at most five elements (two of them are isomorphic, as pointed out to me by Van\v zurov\'a).

In the sequels \cite{IN,KNN,Na,No2,No3,No4,No5,No6,No7,No8}, Nobusawa and his collaborators concentrated mainly on the structure of simple involutory quandles. The exceptions are \cite{KNN} on the structure of latin involutory quandles (see below), \cite{No6} on a close connection between nilpotence and solvability of a quandle and of its displacement group, \cite{No7} on quandles which do not embed into conjugation quandles, and \cite{No8} showing an analogy to the group-theoretical Jordan-H\"older theorem. (In \cite{No4,No6,No7,No8}, Nobusawa already considers general quandles, not necessarily involutory.)

\subsection{Pierce's structure theory}

A big portion of structure theory was developed by Pierce in \cite{Pi1,Pi2}. Both papers are very comprehensive and we will not try to survey the contents here.
The fundamental idea is a correspondence between involutory quandles and groups generated by involutions, via the conjugation representation using the subquandle $Q=\{x\in G:x^2=1\}$ of a group $G$ (i.e. the operation on $Q$ is $x*y=xyx$).
Perhaps most importantly, Pierce obtains an orbit decomposition theorem for involutory quandles in terms of his group representation.
Among other results, he describes free involutory quandles, investigates properties of cycles, and develops a graph theoretical method to study \emph{balanced} involutory quandles (satisfying $xy=y$ iff $yx=x$; for example, cores are balanced).

\subsection{Simple involutory quandles}

In \cite{No2}, Nobusawa points out that the displacement group of a simple involutory quandle $Q$ is a minimal normal
subgroup of $\lmlt Q$, hence, it is either simple, or a direct product of two simple groups which are conjugate. Moreover, in the latter case, $|\dis Q|=|Q|^2$ and $\lmlt Q$ is primitive (primitivity has been shown later in \cite{Na}). The article concludes with an application towards the structure of symmetric groups.

In the next paper \cite{IN}, Ikeda and Nobusawa investigate certain subquandles of the core of the groups $SL_n(q)$ and calculate their displacement groups. They also show that the core of $SL_n(p^k)$ is connected if and only if $p\neq2$ or $n$ is odd; else, it consists of two orbits.

Nobusawa's next paper, \cite{No3}, is concerned with the following construction which often results in a simple involutory quandle.
Let $V$ be a vector space over a field $F$ of characteristic $\neq2$, and let $g$ be a non-degenerate symmetric bilinear form on $V$. 
Let $A=\{a\in V:g(a,a)=1\}$. For $a\in A$, let $L_a$ be the reflection in the hyperplane orthogonal to $a$. 
Let $B$ be the set of all 1-dimensional subspaces of $V$ generated by an element of $A$. Then $Q=(B,*)$ with $\langle a\rangle*\langle b\rangle=L_a(\langle b\rangle)$ is an involutory quandle. It is connected if there exist linearly independent vectors $a,b\in V$ such that $g(a,a)=g(b,b)\neq0$. Then $\lmlt Q$ is primitive if $\dim V>4$ or if $\dim V=4$ and $F\neq\mathbb F_3$.

In \cite{Na}, Nagao discusses involutory quandles with transitive or primitive multiplication groups.
If $\lmlt Q$ is primitive, then $Q$ is simple. If $Q$ is simple, then $\lmlt Q$ is transitive (i.e. $Q$ is connected).
If $G$ is a group generated by a set of involutions $Q$ such that $Q$ is a conjugacy class and the subgroup generated by $QQ=\{ab:\,a,b\in Q\}$ is a minimal normal subgroup of $G$, then the set $Q$ with the operation $x*y=xyx$ is a simple involutory quandle.

In \cite{No4}, Nobusawa extends some of the results on involutory quandles to general quandles, and proves a few new ones. Perhaps the most important fact is that a quandle $Q$ is simple if and only if it is connected, faithful and $\dis Q$ is a minimal normal subgroup of $\lmlt Q$. 
An application towards simplicity of groups (via simplicity of their conjugation quandles) is presented. Continuing in this direction, in \cite{No5} he obtains a new proof that orthogonal groups are simple.

The study of simple quandles culminated by Joyce's classification \cite{Jo2}, including an explicit statement which of them are involutory. The medial simple involutory quandles are the cores of simple abelian groups, and the non-medial ones are realized as a conjugacy class of involutions in a simple non-abelian group. A different description of simple quandles can be found in \cite{AG} (without an explicit reference to involutority).

\subsection{Latin involutory quandles}

There are at least four independent sources of origin for latin involutory quandles. We briefly mention the main ideas and refer to the recent survey \cite{Sta-latin} on latin quandles for a full account.

As noted earlier, they are mentioned in the works of Nobusawa and his followers. In \cite{No1}, latin involutory quandles with prime and prime square elements are proved to be the cores of abelian groups. Most notably, the paper \cite{KNN} by Kano, Nagao and Nobusawa is fully concerned with latin involutory quandles. Using a conjugation representation (as in Pierce's papers), and with the help of Glauberman's $Z^*$-theorem, they prove that a finite involutory quandle $Q$ is latin if and only if the derived subgroup $\lmlt Q'$ has odd order.
They conclude that latin involutory quandles are solvable (informally, they can be constructed by a chain of extensions by medial quasigroups), and satisfy analogies of the Lagrange and Sylow theorems.

Independently, latin involutory quandles appeared in the studies of self-distributive quasigroups by Belousov's group (they called them left symmetric left distributive quasigroups \cite{Bel-book}). In particular, Galkin's theory applies non-trivially to latin involutory quandles \cite[Section 6]{Sta-latin}. Galkin mentions the involutory case explicitly on two occasions.
In \cite{Gal-solv}, he shows that finite latin involutory quandles are solvable, and an infinite latin involutory quandle $Q$ is solvable if and only if for every $\varphi\in\lmlt Q'$ and every $a\in Q$ the mapping $L_a\varphi$ has a unique fixed point.
This is not always the case, as e.g. in $(Q,*)$ with 
\[Q=\{(x_1,x_2,x_3)\in\mathbb{R}^3: x_1^2+x_2^2-x_3^2=-1, x_3>0\}\quad\text{and}\quad x*y=2\langle x,y\rangle x-y\] 
where $\langle x,y\rangle = -x_1y_1-x_2y_2+x_3y_3$.
In \cite{Gal-sub}, Galkin proves that every finite solvable latin quandle has the Lagrange property, but not necessarily the Sylow property. 
In \cite{Gal-sylow}, he proves the Sylow property under the additional assumption that the size of the quandle, and the order of its translations, are coprime; this is always true in the involutory case. 

Another important fact is that there is a one-to-one correspondence between latin involutory quandles and uniquely 2-divisible Bruck loops, see \cite[Setion 5]{Sta-latin} for a precise statement
(Bruck loops are, in a sense, ``the abelian loops among the Bol loops"; a deep theory of uniquely 2-divisible Bruck loops has been developed in \cite{Gl}). 
This connection has a rich history. It was first realized by Robinson in his 1964 PhD thesis, but published only 15 years later in \cite{Rob}. Independently, Belousov and Florya \cite[Theorem 3]{BF} noticed that latin involutory quandles are isotopic to Bol loops, but they did not formulate the full correspondence.
Independently, the correspondence was found by Kikkawa \cite{Ki1} and by Nagy and Strambach \cite{NS}.

Note that commutative involutory quandles are latin. They are precisely the distributive Steiner quasigroups, which correspond to combinatorial designs called \emph{Hall triple systems} \cite{Ben}. It was proved by Bruck \cite{Br} (see also \cite[Lemma 8.5]{Bel-book}) that distributive Steiner quasigroups (i.e. commutative involutory quandles) are precisely the cores of commutative Moufang loops of exponent 3.

\subsection{Medial involutory quandles}

The primary examples of medial involutory quandles are the cores of abelian groups. 
A faithful involutory quandle is medial if and only if it is the core of an abelian group \cite{KNN}, but a plenty of other non-faithful examples exist, see \cite{Ro3} or \cite[Section 7]{JPSZ}.

Medial involutory quandles were studied extensively by Roszkowska-Lech \cite{Ro1,Ro2,Ro3,Ro4} under the name {\it SIE groupoids} (she was probably unaware of the quandle context of her subject). Her main result is that medial involutory quandles can be assembled from pieces that are cores of abelian groups, linking the pieces together by a mesh of group homomorphisms \cite{Ro3}. The construction was recently generalized to all medial quandles \cite{JPSZ}.
She also described the lattice of subvarieties which is isomorphic to the lattice of integers under division with an additional top element \cite{Ro1}, and classified subdirectly irreducible medial involutory quandles \cite{Ro4}.

The subject of \cite{RR} is a special class of reductive medial involutory quandles. The results were also put into a broader context in \cite[Section 6]{JPSZ}.
Another special class, so called group related symmetric groupoids, was studied by Endres in \cite{En1,En3}. 

\subsection{Free involutory quandles}

Free involutory quandles are described in \cite{Pi1}, an explicit normal form of terms is determined in \cite{Sta-gop}. Joyce in \cite{Jo1} shows that the free involutory quandle with two generators is the core of the group of integers, and the free medial involutory quandle with $n$ generators is the subquandle of the core of the group $\Z^{n-1}$ generated by $(0,\dots,0)$, $(1,0,\dots,0),\dots,(0,\dots,0,1)$. It follows that the equational theory of medial involutory quandles is that of the cores abelian groups. Niebrzydowski and Przytycki \cite{NP} studied free objects in the varieties of involutory quandles satisfying $xyx=y$ and $xyxy=y$ (the former identity is commutativity, hence they were partly reproving some older results on Steiner quasigroups).

\subsection{The non-idempotent case}

``Non-idempotent quandles" are referred to as \emph{racks}. The only project specifically addressing non-idempotent involutory racks I know about is the one carried out in our group in Prague more than a decade ago. 
In \cite{Kep}, my former adviser Tom\'a\v s Kepka suggested to study involutory racks (in his terms, left symmetric left distributive groupoids), and stated a few initial observations on their structure and on the lattice of subvarieties. I picked up the proposal and it resulted in some of my earliest papers. In \cite{JKS}, we found a description of non-idempotent subdirectly irreducible involutory racks and obtained some enumeration results. In \cite{Sta-LDLI}, I managed to generalize the subvariety observations into a broader context.
In \cite{Sta-gop}, I discussed several involutory rack operations on groups and derived normal forms for terms in several varieties of involutory racks and quandles. The whole project is summarized in \cite{Sta-A2}.

\section{Connected involutory quandles}

Recall that an involutory quandle is called connected if $\lmlt Q$ is transitive on $Q$. The main result of \cite{HSV} is a canonical correspondence between connected quandles and certain configurations in transitive groups, called quandle envelopes. Fix a set $Q$ and an element $e\in Q$. A quandle envelope is a pair $(G,\zeta)$ such that $G$ is a transitive group on $Q$ and $\zeta\in Z(G_e)$ such that $\langle\zeta^G\rangle=G$. It is easy to see that an envelope $(G,\zeta)$ corresponds to an involutory quandle if and only if $\zeta$ is an involution. 

Consequently, one can easily use all the techniques described in \cite{HSV} to study connected involutory quandles. For example, using their algorithm, one can enumerate connected and latin involutory quandles up to size 47, see Table \ref{Tab:Enum}.

\begin{table}[ht]
$$
\begin{array}{r|rrrrrrrrrrrrrrrr}
n &      1&2&3&4&5&6&7&8&9&10&11&12&13&14&15&16\\\hline
q(n)&    1& 0& 1& 0& 1& 1& 1& 0& 2& 1& 1& 3& 1& 0& 4& 0\\
\ell(n)& 1& 0& 1& 0& 1& 0& 1& 0& 2& 0& 1& 0& 1& 0& {\bf 2}& 0\\
a(n)&    1& 0& 1& 0& 1& 0& 1& 0& 2& 0& 1& 0& 1& 0& 1& 0\\\\
n &      17&18&19&20&21&22&23&24&25&26&27&28&29&30&31&32\\\hline
q(n)&    1& 3& 1& 3& 4& 0& 1& 10& 2& 0& 8& 2& 1& 10& 1& 0\\
\ell(n)& 1& 0& 1& 0&{\bf 2}& 0& 1&  0& 2& 0&{\bf 7}& 0& 1&  0& 1& 0\\
a(n)&    1& 0& 1& 0& 1& 0& 1&  0& 2& 0& 3& 0& 1&  0& 1& 0\\\\
n &      33&34&35&36&37&38&39&40&41&42&43&44&45&46&47& \\\hline
q(n)&    2& 0& 1&16& 1& 0& 2& 8& 1& 8& 1& 0&13& 0& 1\\
\ell(n)& {\bf 2}& 0& 1& 0& 1& 0&{\bf 2}& 0& 1& 0& 1& 0&{\bf 5}& 0& 1\\
a(n)&    1& 0& 1& 0& 1& 0& 1& 0& 1& 0& 1& 0& 2& 0& 1\\
\end{array}
$$
\caption{The numbers $q(n)$ of connected involutory quandles, $\ell(n)$ of latin involutory quandles, and $a(n)$ of connected affine involutory quandles of size $n\le 47$ up to isomorphism.}
\label{Tab:Enum}
\end{table}

\end{document}